\documentclass{agtart_a}
\pdfoutput=1

\usepackage{amscd}

\title[Calabi quasi-morphisms for some non-monotone symplectic manifolds]{Calabi quasi-morphisms for some non-monotone\\symplectic manifolds}
\author{Yaron Ostrover}
\givenname{Yaron}
\surname{Ostrover}
\address{School of Mathematical Sciences\\
Tel Aviv University\\\newline
69978 Tel Aviv\\Israel}
\email{yaronost@post.tau.ac.il}



\received{9 September 2005}
\revised{}
\accepted{15 February 2006}
\proposed{}
\seconded{}
\publishedonline{12 March 2006}
\published{12 March 2006}

\volumenumber{6}
\issuenumber{}
\publicationyear{2006}
\papernumber{14}
\startpage{405}
\endpage{434}

\doi{}
\MR{}
\Zbl{}

\arxivreference{math.SG/0508090}

\AtBeginDocument{\let\bar\wbar\let\tilde\wtilde\let\hat\what}

\makeatletter
\def\cnewtheorem#1[#2]#3{\newtheorem{#1}{#3}[section]
\expandafter\let\csname c@#1\endcsname\c@lemma}

\theoremstyle{plain}
\newtheorem{lemma}{Lemma}[section]
\cnewtheorem{proposition}[lemma]{Proposition}
\cnewtheorem{theorem}[lemma]{Theorem}
\cnewtheorem{corollary}[lemma]{Corollary}

\theoremstyle{remark}
\cnewtheorem{remark}[lemma]{Remark}
\cnewtheorem{definition}[lemma]{Definition}

\makeatother  %  move after \newtheorem block

\newcommand{\Id}{{{\mathchoice {\rm 1\mskip-4mu l} {\rm 1\mskip-4mu l}
      {\rm 1\mskip-4.5mu l} {\rm 1\mskip-5mu l}}}}

\makeatletter \@addtoreset {equation}{section}
\renewcommand\theequation
  {\ifnum \c@subsection>\z@ \arabic{section}.\arabic{subsection}.\arabic{equation}
  \else \arabic{section}.\arabic{equation} \fi}
\makeatother

\begin{document}

\begin{asciiabstract}
In this work we construct Calabi quasi-morphisms on the universal
cover of the group Ham(M) of Hamiltonian diffeomorphisms for some
non-monotone symplectic manifolds. This complements a result by Entov
and Polterovich which applies in the monotone case. Moreover, in
contrast to their work, we show that these quasi-morphisms descend to
non-trivial homomorphisms on the fundamental group of Ham(M).
\end{asciiabstract}

\begin{abstract}
In this work we construct Calabi quasi-morphisms on the universal
cover $\widetilde{\mathrm{Ham}}(M)$ of the group of Hamiltonian
diffeomorphisms for some non-monotone symplectic manifolds. This
complements a result by Entov and Polterovich which applies in the
monotone case. Moreover, in contrast to their work, we show that these
quasi-morphisms descend to non-trivial homomorphisms on the
fundamental group of $\mathrm{Ham}(M)$.
\end{abstract}

\maketitle

\section{Introduction and results}

Let $(M,\omega)$ be a closed connected symplectic manifold of
dimension $2n$. Let\break Ham$(M,\omega)$ denote the group of
Hamiltonian diffeomorphisms of $(M,\omega)$ and let $ {\widetilde
{\rm Ham}(M,\omega)} $ be its universal cover. A celebrated result
by Banyaga \cite{B} states that for a closed symplectic manifold,
Ham$(M,\omega)$ and $ {\widetilde {\rm Ham}(M,\omega)} $ are
simple groups and therefore they do not admit any non-trivial
homomorphism to ${\mathbb R}$. However, in some cases, these
groups admit non-trivial homogeneous quasi-morphisms to ${\mathbb
R}$. Recall that a (real-valued) quasi-morphism of a group $G$ is
a map $r \co G\to {\mathbb R}$ satisfying the homomorphism
equation up to a bounded error, i.e.\ there exists a constant $ C
\geq 0$ such that
$$ |r(g_1g_2) - r(g_1) - r(g_2) | \leq C, \ \ {\rm for \ every \ }
 g_1,g_2 \in G. $$ A quasi-morphism $r$ is called homogeneous if
 $r(g^n) = n r(g)$ for all $g \in G$ and $n \in {\mathbb Z}$. The
 existence of homogeneous quasi-morphisms on the group of Hamiltonian
 diffeomorphisms and/or its universal cover is known for some classes
 of closed symplectic manifolds (see e.g.\ Barge--Ghys~\cite{BarGhy},
 Entov~\cite{E}, Gambaudo--Ghys~\cite{GamGhy} and
 Givental~\cite{Giv}). In a recent work~\cite{EP}, Entov and
 Polterovich showed -- by using Floer and Quantum homology -- that for
 the class of symplectic manifolds which are monotone and whose
 quantum homology algebra is semi-simple, ${\widetilde {\rm
 Ham}(M,\omega)} $ admits a homogeneous quasi-morphism to ${\mathbb
 R}$. In addition to constructing such a quasi-morphism, Entov and
 Polterovich showed that its value on any diffeomorphism supported in
 a Hamiltonianly displaceable open subset equals to the Calabi
 invariant of the diffeomorphism (see Section $2$ below for precise
 definitions). A quasi-morphism with this property is called a Calabi
 quasi-morphism.

The notion ``quasi-morphism" first appeared 
works of Brooks~\cite{Br} and Gromov~\cite{G1} on bounded cohomology
of groups. Since then, quasi-morphisms have become an important tool
in the study of groups. For example, the mere existence of a
homogeneous quasi-morphism on a group $G$ which does not vanish on the
commutator subgroup $G'$ implies that the commutator subgroup has
infinite diameter with respect to the commutator norm (see e.g.\
Bavard~\cite{Bav}). Two well known examples of quasi-morphisms are the
Maslov quasi-morphism on the universal cover of the group of linear
symplectomorphisms of ${\mathbb R}^{2n}$, and the rotation
quasi-morphism defined on the universal cover of the group of
orientation-preserving homeomorphisms of $S^1$. We refer the readers
to Bavard~\cite{Bav} and Kotschick~\cite{K} and the references cited
therein for further details on this subject.  Recently, Biran, Entov
and Polterovich~\cite{BEP}, and Entov and Polterovich~\cite{EP1}
established several other applications of the existence of a Calabi
quasi-morphism regarding rigidity of intersections in symplectic
manifolds. An example of this type is given in
\fullref{intersection-of-stems} and
\fullref{corr-chap-1} below.

In view of the work by Entov and Polterovich \cite{EP}, it is
natural to ask which classes of symplectic manifolds admit a
Calabi quasi-morphism. In a very recent work, Py~\cite{Py}
constructed a homogeneous Calabi quasi-morphism for closed
oriented surfaces with genus greater than $1$. In this note we
concentrate on the case of non-monotone symplectic manifolds. We
will provide some examples of non-monotone rational ruled surfaces
admitting a Calabi quasi-morphism. More precisely, let
 $$ X_{\lambda} = ( S^2 \times S^2, \ \omega_{\lambda} =
\omega \oplus  \lambda \, \omega), \ \ \ 1 \leq \lambda \in
{\mathbb R}, $$ where $\omega$ is the standard area form on the
two-sphere $S^2$ with area $1$, and let
$$ Y_{\mu} = ( {\mathbb C}P^2 \# {\overline {{\mathbb C}P^2}}, \
\omega_{\mu}), \ \ \ 0 < \mu < 1, $$
 be the symplectic blow-up of ${\mathbb C}P^2$ at one point
(see e.g.\ McDuff~\cite{Mc1}, Polterovich~\cite{P3}), where $\omega_{\mu}$ takes the
value $\mu$ on the exceptional divisor, and $1$ on the class of
the line $[{\mathbb C}P^1]$. The manifold $Y_{\mu}$ is the region
$$ \left \{ (z_1,z_2) \in {\mathbb C}^2 \ | \ \mu \leq |z_1|^2 +
|z_2|^2 \leq 1 \right \} $$
with the bounding spheres collapsed
along the Hopf flow. It is known that any symplectic form on these
manifolds is, up to a scaling by a constant, diffeomorphic to one
of the above symplectic forms (see Lalonde--McDuff~\cite{LM3}).

In the monotone case where $\lambda =1$ and $\mu = {\frac 1 3}$,
Entov and Polterovich~\cite{EP} proved the existence of a
homogeneous Calabi quasi-morphism on the universal covers of Ham$
\left (X_{\lambda} \right )$ and Ham$ \left (Y_{\mu} \right )$.
Moreover, they shows that these quasi-morphisms are Lipschitz with
respect to Hofer's metric. For the precise definition of the
Lipschitz property of a quasi-morphism, see Section $2$ below.
Here we prove the following:

\begin{theorem}
\label{Main_Examples}
Let $(M,\omega)$ be one of the following symplectic manifolds:
\renewcommand{\labelenumi}{(\roman{enumi})}
\noindent
\begin{enumerate}
\item $X_{\lambda} = (S^2 \times S^2, \omega_{\lambda})$, \ \
where \ $1 < \lambda \in {\mathbb Q}$.
\item $Y_{\mu} = ( {\mathbb C}P^2 \# {\overline {{\mathbb C}P^2}},  \omega_{\mu})$, \ \ where \
${\frac 1 3} \neq \mu \in  {\mathbb Q} \cap (0,1)$.
\end{enumerate}
\renewcommand{\labelenumi}{\arabic{enumi}.}
Then there exists a homogeneous Calabi quasi-morphism
 $\widetilde r \co{{\widetilde {\rm Ham}}(M,\omega)}\to
{\mathbb R}$, which is Lipschitz with respect to Hofer's metric.
\end{theorem}

It can be shown in the monotone case that any homogeneous
quasi-morphism on the universal cover of Ham$(X_1)$ descends to a
quasi-morphism on Ham$(X_1)$ itself~\cite{EP}. This is due to the
finiteness of the fundamental group $\pi_1 \bigl ( {\rm Ham}(X_1)
\bigr )$, which was proved by Gromov in~\cite{G}. He also pointed
out that the homotopy type of the group of symplectomorphisms of
$S^2 \times S^2$ changes when the spheres have different areas.
McDuff~\cite{Mc}, and Abreu and McDuff~\cite{AM}, showed that the
fundamental groups,  $\pi_1 \bigl ( {\rm Ham}(X_{\lambda}) \bigr
)$ and $ \pi_1 \bigl ( {\rm Ham}(Y_{\mu}) \bigr )$, contain
elements of infinite order for every $0 < \mu < 1$ and for every
$\lambda> 1$. Thus, the above argument will no longer hold in
these cases. Furthermore we claim:

\begin{theorem}
\label{Theorem-Restriction} Let $M$ be one of the manifolds listed
in \fullref{Main_Examples}. Then the restriction of the above
mentioned Calabi quasi-morphism $\widetilde r \co {\widetilde {\rm
Ham}}(M,\omega) \to {\mathbb R}$ to the fundamental group $\pi_1
\bigl ( {\rm Ham}(M,\omega) \bigr ) \subset {\widetilde {\rm
Ham}}(M,\omega)$ gives rise to a non-trivial homomorphism.
\end{theorem}

This differs from the situation described in~\cite{EP} where it
was proven that for $M = {\mathbb C}P^n$ endowed with the
Fubini--Study form, or for $M = S^2 \times S^2$ equipped with the
split symplectic form $\omega \oplus \omega$, the restriction of
the Calabi quasi-morphism to the fundamental group $\pi_1 \bigl (
{\rm Ham}(M) \bigr )$ vanishes identically.

For technical reasons, we shall assume in what follows that $M$ is
a rational strongly semi-positive symplectic manifold. Recall that
a symplectic manifold $M$ is rational if the set $\{ \omega(A) \ |
\ A \in \pi_2(M)\}$ is a discrete subset of ${\mathbb R}$ and
strongly semi-positive if for every $A \in \pi_2(M)$ one has
$$ 2-n \leq c_1(A) < 0 \Longrightarrow \omega(A) \leq 0,$$
where $c_1 \in H^2(M,{\mathbb Z})$ denotes the first Chern class
of $M$. The assumption that $M$ is strongly semi-positive is a
standard technical assumption (see e.g.\ Piunikhin--Salamon--Schwarz~\cite{PSS},
Seidel~\cite{Se})
which guarantees, roughly speaking, the good-behavior of some
moduli spaces of pseudo-holomorphic curves. Note that every
symplectic manifold of dimension $4$ or less, in particular the
manifolds listed in \fullref{Main_Examples}, is strongly
semi-positive. The rationality assumption is also a technical
assumption. It plays a role, for example, in \fullref{PD_lemma}
below, where for non-rational symplectic manifolds the action
spectrum is a non-discrete subset of ${\mathbb R}$ and our method
of proof fails.

In fact, the examples in \fullref{Main_Examples} are special
cases of a more general criterion for the existence of a Calabi
quasi-morphism. In~\cite{EP}, such a criterion was given for
closed monotone symplectic manifolds. This criterion is based on
some algebraic properties of the quantum homology algebra of
$(M,\omega)$. More precisely, recall that as a module the quantum
homology of $M$ is defined as $QH_*(M) = H_*(M) \otimes \Lambda$,
where $\Lambda$ is the standard Novikov ring $$ \Lambda = \Bigl \{
\sum_{A \in \Gamma} \lambda_A q^A \ | \
 \lambda_A \in {\mathbb Q}, \  \#  \{ A \in \Gamma \ | \ \lambda_A \neq 0, \ \omega(A) > c \}
 < \infty, \ \forall \, c \in {\mathbb R} \Bigr \}.
$$
Here $\Gamma = \pi_2(M) ~ / ~ ({\rm ker} \, c_1 \cap {\rm ker} \,
\omega )$, where $c_1$ is the first Chern class. A grading on
$\Lambda$ is given by ${\rm deg} (q^A) = 2c_1(A)$. We shall denote
by $\Lambda_k$ all the elements in $\Lambda$ with degree $k$. We
refer the readers to McDuff--Salamon~\cite{MS2} and to Subsection $3.1$ below for
a more detailed exposition and for the precise definition of the
quantum product on $QH_*(M)$. In the monotone case, i.e.\ where
there exists $\kappa
> 0$ such that $\omega = \kappa \cdot c_1$, 
the Novikov ring $\Lambda$ can be identified with the field of
Laurent series $\sum \alpha_j x^j$, with coefficients in ${\mathbb
Q}$, and all $\alpha_j$ vanish for $j$ greater than some large
enough $j_0$. 
In this case we say that the quantum homology $QH_*(M)$ is {\it
semi-simple} if it splits with respect to multiplication into a
direct sum of fields, all of which are 
finite dimension linear spaces over $\Lambda$. It was shown
in~\cite{EP} that for monotone symplectic manifolds with
semi-simple quantum homology algebra there exists a Lipschitz
homogenous Calabi quasi-morphism on the universal cover of the
group of Hamiltonian diffeomorphisms.

In the non-monotone case the above definition of semi-simplicity
will no longer hold since $\Lambda$ is no longer a field. However,
it turns out that a similar criterion to the above still exists in
this case. More precisely, we focus upon the sub-algebra
$QH_{2n}(M) \subset QH_*(M)$ over the sub-ring $\Lambda_0 \subset
\Lambda$. This sub-algebra is the degree component of the identity
in $QH_*(M)$. Using the fact that in the non-monotone case the
sub-ring $\Lambda_0$ can be identified with the field of Laurent
series, we say as before that $QH_{2n}(M)$ is semi-simple over
$\Lambda_0$ if it splits into a direct sum of fields with respect
to multiplication. Denote by $N_M$ the minimal Chern number of $M$
defined as the positive generator of the image $ c_1 \left (
\pi_2(M) \right ) \subseteq {\mathbb Z}$ of the first Chern class
$c_1$. The following criterion is a generalization of Theorem
$1.5$ from~\cite{EP} to the rational strongly semi-positive case.

\begin{theorem}
\label{The criterion}
Let $(M,\omega)$ be a closed connected rational strongly
semi-positive symplectic manifold of dimension $2n$. Suppose that
the quantum homology subalgebra $QH_{2n}(M) \subset QH_*(M)$ is a
semi-simple algebra over the field $\Lambda_0$ and that 
$N_M$ divides $n$.  Then there exists a Lipschitz homogeneous
Calabi quasi-morphism $\widetilde r \co {{\widetilde {\rm
Ham}}(M,\omega)} \to {\mathbb R}.$
\end{theorem}

For the manifolds $X_{\lambda}$ and $Y_{\mu}$ listed in
\fullref{Main_Examples} the minimal Chern number $N_M$ is $2$
and $1$ respectively. Thus, one of our main tasks is to prove that
for these manifolds the top-dimension quantum homology subalgebra
$QH_{4}(M)$ is semi-simple over the field $\Lambda_0$.

As a by-product of \fullref{Main_Examples}, we generalize a
result regarding rigidity of intersections obtained by Entov and
Polterovich in~\cite{EP1}. To describe the result, we recall first
the following definitions. For a symplectic manifold $M$ denote by
$\{ \cdot , \cdot \}$ the standard Poisson brackets on
$C^{\infty}(M)$. A linear subspace ${\cal A} \subset
C^{\infty}(M)$ is said to be Poisson-commutative if $\{F,G\}=0$
for all $F,G \in {\cal A}$. We associate to a finite-dimensional
Poisson-commutative subspace ${\cal A} \subset C^{\infty}(M)$ its
moment map $\Phi_{\cal A} \co M \to {\cal A}^*$, defined by $
\langle \Phi_{\cal A}(x),F \rangle = F(x)$. A non-empty subset of
the form $\Phi_{\cal A}^{-1}(p), \ p \in {\cal A}^*$, is called a
fiber of ${\cal A}$. A fiber $X \subset M$ is said to be
displaceable if there exists a Hamiltonian diffeomorphism $\varphi
\in {\rm Ham}(M)$ such that $\varphi(X) \cap X = \emptyset$. The
following definition was introduced in~\cite{EP1}:

\begin{definition}
A closed subset $X \subset M$ is called a {\it stem}, if there
exists a finite-dimensional Poisson-commutative subspace ${\cal A}
\subset C^{\infty}(M)$, such that $X$ is a fiber of ${\cal A}$ and
each fiber of ${\cal A}$, other than $X$, is displaceable.
\end{definition}

In Theorem $2.4$ of~\cite{EP1}, Entov and Polterovich showed that
for a certain class of symplectic manifolds, any two stems have a
non-empty intersection. What they used, in fact, was only the
existence of a Lipschitz homogeneous Calabi quasi-morphism for
manifolds in this class. Using the exact same line of proof, the
following theorem follows from \fullref{Main_Examples} above.

\begin{theorem} \label{intersection-of-stems}
Let $M$ be one of the manifolds listed in
\fullref{Main_Examples}. Then any two stems in $M$ intersect.
\end{theorem}

An example of a stem in the case where $M=X_{\lambda}$ is the
product of two equators. More precisely, 
we identify $X_{\lambda}$ with ${\mathbb C}P^1 \times {\mathbb
C}P^1$ in the obvious way. Denote by $L \subset X_{\lambda}$ the
Lagrangian torus defined by
$$ L = \left \{ \, ( [z_0:z_1],[w_0:w_1] ) \in {\mathbb C}P^1 \times
{\mathbb C}P^1  \mid \ \ |z_0| = |z_1|, \ |w_0| = |w_1| \right \}
$$ The proof that $L$ is a stem goes along the same line as Corollary $2.5$
of~\cite{EP1}. Since the image of a stem under any
symplectomorphism of $M$ is again a stem we get:
\begin{corollary} \label{corr-chap-1}
Let $X_{\lambda}$ be one of the manifolds in the first class of
manifolds listed in \fullref{Main_Examples} above. Then for
any symplectomorphism $\varphi$ of $X_{\lambda}$ we have $L \cap
\varphi(L) \neq 0$.
\end{corollary}

{\bf Organization of the paper}\qua In Section $2$ we
recall some definitions and notations related to the Calabi
quasi-morphism. In Section $3$ we briefly review the definition of
the quantum homology algebra $QH_*(M)$. We then describe the
quantum homologies of our main examples and state some of their
properties. In Section $4$ we recall the definition of Floer
homology and some relevant notions. Section $5$ is devoted to the
proof of \fullref{Main_Examples} and \fullref{The
criterion}. In Section $6$ we discuss the restriction of the
Calabi quasi-morphisms on the fundamental group of ${\rm Ham}(M)$.
In Section $7$ we prove \fullref{Theorem-Restriction} and in
the last section we prove the Poincar\'{e} duality type lemma
which is stated and applied in Section $5$.

{\bf Acknowledgments}\qua I would like to thank my supervisor,
Professor Leonid Polterovich, for his guidance, patience and
invaluable advice during the preparation of this work. I would also
like to thank Professor Michael Entov for many helpful suggestions. I
am grateful to Professor Dusa McDuff for many useful remarks, and for
inviting me to Stony Brook University and giving me the opportunity to
present this work, as well as for her warm hospitality. A part of this
work was done during my visit to the Max-Planck Institute in
Leipzig. I thank the Max-Planck Institute for their invitation and
especially I would like to thank Professor Matthias Schwarz for
enlightening conversations and his warm hospitality.  This research
was partially supported by the Israel Science Foundation grant \#
11/03.

\section{Preliminaries on Calabi quasi-morphism}

\setcounter{equation}{0}

In this section we recall the definition of a Calabi
quasi-morphism introduced in~\cite{EP}. We start with the
definition of the classical Calabi invariant (see Banyaga \cite{B} and
Calabi \cite{C}). Let $(M,\omega)$ be a closed connected symplectic
manifold. Given a Hamiltonian function $H \co S^1 \times M
\rightarrow {\mathbb R}$, set $H_t := H(t,\cdot)$ and denote by
$\varphi$ the time-1-map of the Hamiltonian flow
$\{\varphi_H^t\}$. The group of Hamiltonian diffeomorphisms ${\rm
Ham}(M,\omega)$ consists of all such time-1-maps. 
Let $ {\widetilde {\rm Ham}(M,\omega)} $ be the universal cover of
Ham$(M,\omega)$. For a non-empty open subset $U$ of $M$, we denote
by $ \widetilde {\rm Ham}_U(M,\omega)$  the subgroup of $
\widetilde {\rm Ham}(M,\omega)$, consisting of all elements that
can be represented by a path $\{ \varphi_{H}^t \}_{t
\in [0,1]}$ starting at the identity 
and generated by a Hamiltonian function $H_t$ supported in $U$ for
all $t$. For $\varphi \in \widetilde {\rm Ham}_U(M,\omega)$ we
define ${\rm Cal}_U \co \widetilde {\rm Ham}_U(M,\omega) \to
{\mathbb R}$ by
\begin{equation*} \label{Def_Calabi_invariant} \varphi \mapsto
\int_0^1dt \int_M H_t \, \omega^n.\end{equation*}
This map is well defined, i.e.\ it is independent of the specific
choice of the Hamiltonian function generating $\varphi$. Moreover,
it is a group homomorphism called the Calabi homomorphism.

Recall that a non-empty subset $U$ of $M$ is called Hamiltonianly
displaceable if there exists a Hamiltonian diffeomorphism $\varphi
\in {\rm Ham} (M,\omega)$ such that
$ \varphi(U) \cap {\rm Closure}(U) = \emptyset.$ The following two
definitions were introduced in~\cite{EP}.

\begin{definition} \label{Calabi_prop}
A quasi-morphism on $\widetilde {\rm Ham}(M,\omega)$ coinciding
with the Calabi homomorphism ${\rm Cal}_U \co \widetilde {\rm Ham}_U
(M,\omega) \rightarrow {\mathbb R}$ on any open and Hamiltonianly
displaceable
set $U$ 
is called a {\it Calabi quasi-morphism}.
\end{definition}

\begin{definition} \label{Leipschitz prop}
A quasi-morphism $ r \co \widetilde {\rm Ham}(M,\omega) \to
{\mathbb R}$ is said to be {\it Lipschitz with respect to Hofer's
metric} if there exists a constant $K > 0$ so that
$$ |r(\varphi_H) - r(\varphi_F)| \leq K \cdot \| H - F
\|_{C^0}$$
\end{definition}

For the relation of $\|H - F \|_{C^0}$ to the Hofer distance
between the corresponding Hamiltonian diffeomorphisms $\varphi_H$
and $\varphi_H$ see e.g.~\cite{EP}.

\section{The Quantum homology of our main examples}

\subsection{The quantum homology algebra} \label{Sub-Sec_QH}
\setcounter{equation}{0}

In this section we briefly recall the definition of the quantum
homology ring of $(M^{2n},\omega)$. We refer the readers to
\cite{MS2} for a detailed exposition on this subject. Let $M$ be a
closed rational strongly semi-positive symplectic manifold of
dimension $2n$. By abuse of notation, we shall write $\omega(A)$
and $c_1(A)$ for the integrals of $\omega$ and $c_1$ over $A \in
\pi_2(M) $. Let $\Gamma$ be the abelian group
\begin{equation} \label{The_group_Gamma}
\Gamma = \pi_2(M) ~ / ~ ({\rm ker} \, c_1 \cap {\rm ker} \, \omega
).
\end{equation}
We denote by $\Lambda$ the Novikov ring
\begin{equation} \label{Def_Novikov_Ring}
\Lambda = \Bigl \{ \sum_{A \in \Gamma} \lambda_A q^A \ | \
 \lambda_A \in {\mathbb Q}, \  \#  \{ A \in \Gamma \ | \ \lambda_A \neq 0, \ \omega(A) > c \}
 < \infty, \ \forall \, c \in {\mathbb R}\, \Bigr \}.
\end{equation}
This ring comes with a natural grading defined by ${\rm deg}(q^A)
= 2c_1(A)$. We shall denote by $\Lambda_k$ all the elements of
$\Lambda$ with degree $k$. Note that  $\Lambda_k = \emptyset$ if
$k$ is not an integer multiple of $2N_M$, where $N_M$ is the
minimal Chern number of $M$ defined by $ c_1 \left ( \pi_2(M)
\right ) = N_M {\mathbb Z}$.

As a module, the quantum homology ring of $(M,\omega)$ is defined
as
$$  QH_{*}(M) = QH_{*}(M,\Lambda) = H_{*}(M,{\mathbb Q}) \otimes \Lambda. $$
A grading on $QH_*(M)$ is given by $ {\rm deg} (a \otimes q^A ) =
{\rm deg} (a) + 2c_1(A),$ where ${\rm deg}(a)$ is the standard
degree of the class $a$ in the singular homology of $M$. Next, we
define the quantum product on $QH_{*}(M)$ (cf~\cite{MS2},
\cite{RT}). For $a \in H_i(M)$ and $b \in H_j(M)$ we define $a * b
\in QH_{i+j-2n}(M)$ as
$$ a * b = \sum_{A \in \Gamma}(a * b)_A \otimes q^{-A}, $$
where $(a * b)_A \in H_{i+j-2n+2c_1(A)}(M)$ is determined by the
requirement that
$$ (a * b)_A \ \circ \ c = \Phi_A(a,b,c) \ \ \ \ \ {\rm for \ all}
\ c \in H_*(M). $$ Here $\circ$ is the usual intersection product
on $H_*(M)$, and $\Phi_A(a,b,c)$ denotes the Gromov--Witten
invariant that counts the number of pseudo-holomorphic curves
representing the class $A$ and intersecting with a generic
representative of each of $a,b,c \in H_*(M)$. The product $*$ is
extended to $QH_{*}(M)$ by linearity over the ring $\Lambda$. Note
that the fundamental class $[M]$ is the unity with respect to the
quantum multiplication.

It follows from the definitions that the zero-degree component of
$a *b$ coincides with the classical cap-product  $a \cap b$ in the
singular homology. Moreover, there exists a natural pairing
$\Delta \co QH_k(M) \times QH_{2n-k}(M) \to \Lambda_0$ defined by
$$ \Delta \left ( \sum a_A \otimes q^{A} ,  \sum b_B \otimes q^{B} \right ) =
\sum_{c_1(A)=0} \Bigl ( \sum_{B} (a_{-B} \circ b_{B+A}) \Bigr )
q^{A}. $$ The fact that the inner sums on the right hand side are
always finite follows from the finiteness condition
in~$\ref{Def_Novikov_Ring}$. Moreover, the pairing $\Delta$
defines a Frobenius algebra structure, i.e.\ it is non-degenerate
in the sense that $\Delta(\alpha,\beta) = 0$ for all $\beta$
implies $\alpha =0$, and $\Delta(\alpha,\beta) = \Delta(\alpha *
\beta, [M])$. Notice that $\Delta$ associates to each pair of
quantum homology classes $\alpha ,\beta  \in QH_*(M)$ the
coefficient of the class $P=[{\rm point}]$ in their quantum product. We
also define a non-degenerate ${\mathbb Q}$--valued pairing $\Pi$ to
be the zero order term of $\Delta$, i.e.\
\begin{equation} \label{Def_Poincare_pairing}
\Pi \left ( \sum a_A \otimes q^{A} ,  \sum b_B \otimes q^{B}
\right ) = \sum_{B} (a_{-B} \circ b_{B}).
\end{equation}
Note that $\Pi(\alpha,\beta) = \Pi(\alpha * \beta, [M])$ for every
pair of quantum homology classes $\alpha$ and $\beta$.
Furthermore, the finiteness condition in the definition of the
Novikov ring $\ref{Def_Novikov_Ring}$ leads to a natural
valuation function $val\co QH_*(M) \to {\mathbb R}$ defined by
\begin{equation} \label{Def_valuation} val \bigl ( \sum_{A \in \Gamma} a_A \otimes q^A  \bigr )
= {\rm max} \{ \omega(A) \ | \ a_A \neq 0 \},  \ \  {\rm and} \ \
val(0) = - \infty. \end{equation}

\subsection{The case of ${S^2 \times S^2}$}
\setcounter{equation}{0}

Let $X_{\lambda} = S^2 \times S^2$ be equipped with the split
symplectic form $\omega_{\lambda}= \omega  \oplus \lambda \,
\omega$, where $\lambda > 1$. In this subsection we discuss
several issues regarding the quantum homology of the manifold
$X_{\lambda}$ and in particular we show that the quantum homology
subalgebra $QH_{4}(X_{\lambda}) \subset QH_*(X_{\lambda})$ is a
field for every $\lambda > 1$.

Denote the standard basis of $H_*(X_{\lambda})$ by $ P = [ {\rm
point}], \  A = [S^2 \times {\rm point}] , \ B = [{\rm point}
\times S^2]$ and the fundamental class $M=[X_{\lambda}]$. The
quantum homology of $X_{\lambda}$ is generated over the Novikov
ring $\Lambda$ by these elements. Since $\lambda
> 1$, it follows that $\Gamma = \pi_2(X_{\lambda})$, where the last is isomorphic to the free abelian group
generated by $A$ and $B$. From the following Gromov--Witten
invariants (see e.g.~\cite{EP},~\cite{MS2}):
$$ \Phi_{A+B} (  P, P, P ) =1, \ \    \Phi_{0} ( A, B, M ) =1, \
\  \Phi_{A} (  P, B, B ) =1, \  \    \Phi_{B} ( P, A, A ) =1,  $$
one finds the quantum identities: \begin{equation}
\label{quantum-relations-of-the-first-example} A
* B = P, \ \ \ A^2 = M \otimes q^{-B},   \ \ \  B^2 = M \otimes
q^{-A}.
\end{equation} Next, instead of the standard basis $\{A,B \}$ of $\Gamma$,
we consider the basis $\{e_1 , e_2 \} = \{B-A,A\}$. Set $x =
q^{e_1}$ and $y=q^{e_2}$. In this notation, the quantum product of
the generators of $QH_*(X_{\lambda})$ becomes
\begin{equation} \label{quantum-relations-of-the-first-example-2}
A * B = P, \ \ \ A^2 = M \otimes x^{-1} y^{-1},   \ \ \  B^2 = M
\otimes y^{-1}. \end{equation}
It follows from the definition of the Novikov
ring~$\ref{Def_Novikov_Ring}$ that
$$ \Lambda =   \left \{
        \sum \lambda_{\alpha, \beta}  \cdot x^{\alpha} y^{\beta} \ \big | \
\lambda_{\alpha, \beta} \in {\mathbb Q}\, \right \},$$ where each
sum satisfies the following finiteness condition:
$$
 \# \left \{
(\alpha,\beta)  \ | \ \lambda_{\alpha, \beta} \neq 0, \
\alpha(\lambda -1) + \beta
> c  \right \} < \infty,  \
 \forall \, c \in {\mathbb R}. $$
Taking into account the above mentioned grading of $\Lambda$ we
get
\begin{eqnarray*}
\Lambda_{4k} & = & \left \{
 \sum \lambda_{\alpha, \beta}  \cdot  x^{\alpha}  y^{\beta}  \ \in \Lambda \ \big | \
 4 \beta =   2c_1(\alpha e_1 + \beta e_2)   = 4k \right \} \\ \\
 & = &  \left \{ \sum \lambda_{\alpha} \cdot x^{\alpha} y^{k} \ \big | \
  {  \  \#  \{ \alpha \   | \ \lambda_{\alpha} \neq 0, \ \alpha(\lambda -1)  > d \}
 < \infty, \ \forall \, d \in {\mathbb R} } \right  \} .
\end{eqnarray*}
The finiteness condition above implies that $\lambda_{\alpha}$
vanishes for large enough $\alpha$'s.

\begin{lemma}  \label{QH_4_is_a_subfield_1}
For any $\lambda > 1$, the subalgebra $QH_4(X_{\lambda}) \subset
QH_*(X_{\lambda})$ is a field.
\end{lemma}

\begin{proof}
Let $0 \neq \gamma \in QH_{4}(X_{\lambda})$.
 Since
$QH_4(X_{\lambda})= H_4(X_{\lambda}) \otimes \Lambda_{0} +
H_0(X_{\lambda}) \otimes \Lambda_{4}$, it follows that
$$ \gamma  =  M \otimes \sum \lambda_{\alpha_1} \ x^{  \alpha_1 } +
P \otimes y \sum \lambda_{\alpha_2} \ x^{ \alpha_2}, $$ where
$\lambda_{\alpha_1}$ and $\lambda_{\alpha_2}$ vanish for large
enough $\alpha_1$ and $\alpha_2$ respectively. Next, let $\beta =
P \otimes y$ be a formal variable. From the above multiplicative
relations~$\ref{quantum-relations-of-the-first-example-2}$, we
see that $\beta^2 = M \otimes x^{-1}$. Hence, we can consider the
following ring identification:
$$
QH_4(X_{\lambda}) \simeq {\cal R}[\beta] ~ / ~ {\cal I}, $$ where
${\cal I}$ is the ideal generated by  $\beta^2 - x^{-1}$ and
${\cal R}= {\mathbb Q}[[x]$ is the ring of Laurent series $\sum
\alpha_j x^j$, with coefficients in ${\mathbb Q}$, and all
$\alpha_j$ vanish for $j$ greater than some large enough $j_0$.
Note that for any Laurent series $\Phi(x) \in {\cal R}$, the
maximal degree of $\Phi^2(x)$ is either zero or even. Therefore
${\cal R}$ does not contain a square root of $x^{-1}$ and hence
$\cal I$ is a maximal ideal. Thus, we conclude that
$QH_4(X_{\lambda})$ is a two-dimensional extension field of ${\cal
R}$. This completes the proof of the lemma.
\end{proof}

\begin{remark}
{\rm Note that the above statement no longer holds in the monotone
case where $\lambda = 1$, since $QH_{4}(X_1)$ contains zero
divisors (see e.g.~\cite{EP},~\cite{MS2}). }
\end{remark}

\subsection{The case of ${\mathbb C}P^2 \# \overline{{\mathbb C}P^2}$}
\setcounter{equation}{0}

Here we study the quantum homology algebra of $Y_{\mu} = (
{\mathbb C}P^2 \# {\overline {{\mathbb C}P^2}}, \omega_{\mu} )$,
which is the symplectic one-point blow-up of ${\mathbb C}P^2$
introduced in Section $1$. We will show that the quantum homology
subalgebra $QH_4(Y_{\mu})$, which plays a central role in the
proof of \fullref{Main_Examples}, is semi-simple. It is worth
mentioning (see \fullref{algebraic_structure_of_QH(Y)} below)
that the algebraic structure of $QH_4(Y_{\mu})$ turns out to be
dependent on $\mu$.

We denote by $E$ the exceptional divisor and by $L$ the class of
the line $[{\mathbb C}P^1]$. Recall that for $0 < \mu < 1$,
$\omega_{\mu}$ is a symplectic form on $Y_{\mu}$ with
$\omega_{\mu}(E) = \mu$ and $\omega_{\mu}(L) = 1$. Denote the
class of a point by $P=[{\rm point}]$ and set $F = L - E$. The elements
$P,E,F$ and the fundamental class $M=[Y_{\mu}]$ form a basis of
$H_*(Y_{\mu})$.

The following description of the multiplicative relations for the
generators of $QH_*(Y_{\lambda})$ can be found in~\cite{Mc1}.

\begin{center}
\begin{tabular}{ll} \label{quantum-relations-blow-up}
$ P * P  = (E+F) \otimes q^{-E-F} $, &  \quad $E * P  =  F \otimes q^{-F}$, \\
$ P * F  =   M \otimes q^{-E-F}$, & \quad $E * E  =  -P + E
\otimes q^{-E} +
M \otimes q^{-F}$, \\
$E * F  =  P-E \otimes q^{-E},  $ & \quad $F * F = E  \otimes
  q^{-E}$. \end{tabular}
\end{center}

Consider the rational non-monotone case where ${\frac 1 3} \neq
\mu \in {\mathbb Q} \cap (0,1)$. Note that in this case $\Gamma
\simeq {\mathbb Z} \otimes {\mathbb Z}$. As in the previous
example of $S^2 \times S^2$, we apply a unimodular change of
coordinates and consider the following basis of $\Gamma$
$$ \Gamma \simeq \left\{
 \begin{array}{l}
 {\rm Span}_{\mathbb Z} \{ F-2E , \ E \}, \quad  0 < \mu < {\frac 1 3}  \\ \\
 {\rm Span}_{\mathbb Z} \{ 2E-F, \ E \},   \quad {\frac 1 3} < \mu < 1
\end{array} \right.
$$
Denote $e_1 = F - 2E$, $e_2 = E$ when $0 < \mu < {\frac 1 3}$ and
 $e_1 = 2E - F$, $e_2 = E$ when ${\frac 1 3} < \mu < 1$. Set $x=q^{e_1}$ and
 $y=q^{e_2}$.
From the  definition of the Novikov
ring~$\ref{Def_Novikov_Ring}$ we have
$$ \Lambda =  \left \{
        \sum
   \lambda_{\alpha, \beta}  x^{\alpha} y^{\beta} \ \big | \
\lambda_{\alpha, \beta} \in {\mathbb Q} \right \},$$ where each
sum satisfies the following finiteness condition:
$$ \# \left \{
(\alpha,\beta)  \ | \ \lambda_{\alpha, \beta} \neq 0, \ \alpha | 3
\mu - 1 | + \beta \mu
> c  \right \} < \infty,  \
 \forall \, c \in {\mathbb R}. $$
The graded Novikov ring has the form
\begin{eqnarray*}
 \Lambda_{2i} & = & \left \{  
 \sum \lambda_{\alpha, \beta}  \cdot  x^{\alpha}  y^{\beta}  \ \in \Lambda \ \big | \
  2 \beta =  2c_1(\alpha e_1 + \beta e_2) =2i \right \} \\ \\
  & = &  \left \{ \sum \lambda_{\alpha} \cdot x^{\alpha} y^{i} \  \big | \
  {  \  \#  \{ \alpha \ | \ \lambda_{\alpha} \neq 0, \ \alpha |3 \mu - 1|  > d \}
 < \infty, \ \forall \, d \in {\mathbb R} } \right  \}.
\end{eqnarray*}
Next we present the quantum product of $QH_*(Y_{\mu})$ with
respect to the above basis of $\Gamma$.

\begin{center}
\begin{tabular}{ll}
$ P * P  = (E+F) \otimes x^{\kappa} y^{-3} $, &  \quad $E * P  =  F \otimes x^\kappa y^{-2}$, \\
$ P * F  =   M \otimes x^\kappa y^{-3}$, & \quad $E * E  =  -P + E
\otimes y^{-1} +
M \otimes x^\kappa y^{-2}$, \\
$E * F  =  P-E \otimes y^{-1},  $ & \quad $F * F = E  \otimes
  y^{-1}$, \end{tabular}
\end{center}
 where $\kappa = {\rm sgn}(3\mu-1)$, i.e.\ $\kappa =1$ for ${\frac
1 3} < \mu < 1$, and $\kappa = -1$ for $0 < \mu  < {\frac 1 3}$.

\begin{lemma}  \label{QH_4_is_a_subfield_2}
The subalgebra $QH_4(Y_{\mu}) \subset QH_*(Y_{\mu})$ is
semi-simple.
\begin{proof}
Since $QH_4(Y_{\mu})=
H_4(Y_{\mu}) \otimes \Lambda_{0} + H_2(Y_{\mu}) \otimes
\Lambda_{2} + H_0(Y_{\mu}) \otimes \Lambda_{4}$, it follows that
for every $0 \neq \delta \in QH_{4}(Y_{\mu})$
\begin{eqnarray*} \delta  & = & M \otimes \sum \lambda_{\alpha_1}
 x^{\alpha_1} + E \otimes  y  \sum \lambda_{\alpha_2} x^{\alpha_2} \\
 & + & F \otimes y \sum \lambda_{\alpha_3}  x^{\alpha_3} + P \otimes
 y^{2} \sum \lambda_{\alpha_4}x^{\alpha_4}, \end{eqnarray*}
where $\lambda_{\alpha_i}$ vanish for large enough $\alpha_i$ for
$i=1,2,3,4$. Next, put $\beta_1 = E \otimes  y$, $\beta_2 = F
\otimes y$ and $\beta_3 = P \otimes y^{2}$. From the above
multiplication table, we see that
\begin{equation*} \left\{
\begin{array}{l}
 \beta_1^2   =  - \beta_3 + \beta_1 + x^\kappa, \\  \beta_2 ^ 2 = \beta_1, \\ \beta_3^2 = x^\kappa(\beta_1 + \beta_2) \\
\beta_1  \cdot  \beta_2  =  \beta_3 - \beta_1, \\ \beta_2 \cdot
\beta_3 =x^\kappa, \\ \beta_1 \cdot \beta_3 = x^\kappa \beta_2.
\end{array} \right.
\end{equation*}
Thus, we have the following ring identification:
$$ QH_4(Y_{\mu}) \simeq {\cal R}[\beta_1, \beta_2, \beta_3]  ~ / ~
{\cal I}, $$ where ${\cal R} = {\mathbb Q}[[x]$ is the ring of
Laurent series $\sum \alpha_j  x^j$ 
and ${\cal I}$ is the ideal generated by the above relations. It
is easy to check that the sixth equation follows immediately from
the second and the fifth equations and hence, it can be
eliminated. Moreover, by isolating $\beta_3$ and $\beta_1$ from
the first and the second equations respectively, we conclude that
the above system is equivalent to the following one:
\begin{equation} \label{mul-rel} \left\{
\begin{array}{l}
  (\beta_2^2 - \beta_2^4 + x^\kappa)^2 = x^\kappa(\beta_2^2 + \beta_2) \\
\beta_2^3  =  - \beta_2^4 + x^\kappa, \\ \beta_2 \cdot (\beta_2^2
- \beta_2^4 + x^\kappa) =x^\kappa,
\end{array} \right.
\end{equation}
Moreover, we claim that in fact 
\[
QH_4(Y_{\mu}) ~ \simeq ~ {\cal R}[\beta_1, \beta_2, \beta_3] ~ / ~
{\cal I} ~ \simeq ~ {\cal R}[\beta_2] ~ / ~  {\cal J} \] where
${\cal J}$ is the ideal generated by $\beta_2^4 + \beta_2^3 -
x^\kappa$. Indeed, the first equation in~$\ref{mul-rel}$ is
obtained by multiplying the third equation by $\beta_2^2+\beta_2$
and assigning the second equation. The third equation is obtained
from the second after multiplying it by $\beta_2-1$. Next, note
that 
the polynomial $\beta_2^4 + \beta_2^3 - x^\kappa$ does not share a
common root with its derivative since the roots of the derivative
are $0$ and $-3/4$. Thus, 
it
has no multiple roots in ${\cal R}$ and hence the quantum homology
subalgebra $QH_4(Y_{\mu})$ is semi-simple as required.
\end{proof}
\end{lemma}

\begin{remark} \label{algebraic_structure_of_QH(Y)} {\rm
Strangely enough, it follows from the above lemma that the
algebraic structure of the quantum homology subalgebra
$QH_4(Y_{\mu})$ depends on $\mu$. More precisely, it can be shown
 that the polynomial $\beta_2^4 + \beta_2^3 - x^\kappa$ is
irreducible over ${\cal R}$ for $\kappa = 1$ while reducible for
$\kappa = -1$. Thus, $QH_4(Y_{\mu})$ is a field when ${\frac 1 3}
< \mu < 1$, while for $ 0 < \mu < {\frac 1 3}$, it is a direct sum
of fields. We omit here the technical details because for our
purpose, it is sufficient that $QH_4(Y_{\mu})$ is semi-simple.  }
\end{remark}

\section{Preliminaries on Floer homology} \label{Pre_on_Floer}

In this section we give a brief review of Floer homology. In
particular we present some definitions and notions which will be
relevant for the proof of our main results. We refer the readers
to~\cite{Sa} or~\cite{MS2} for a more detailed description.

Let $(M,\omega)$ be a closed, connected and strongly semi-positive
symplectic manifold. Let $J = \{ J_t \} _{0 \leq t \leq 1}$ be a
periodic family of $\omega$--compatible almost complex structures.
We denote by ${\cal L}$ the space of all smooth contractible loops
$x \co S^1 = {\mathbb R} / {\mathbb Z} \to M$. Consider a covering
${\widetilde {\cal L}}$ of ${\cal L}$ whose elements are
equivalence classes $[x,u]$ of pairs $(x,u)$, where $x \in {\cal
L}$, $u$ is a disk spanning $x$ in $M$, and where
$$ (x_1,u_1) \sim (x_2,u_2) \ \ \ {\rm if \ and \ only \ if} \  x_1=x_2 \
\ {\rm and} \ \ \omega(u_1 \# u_2) = c_1(u_1 \# u_2) = 0. $$ The
group of deck transformations of $\widetilde {\cal L}$ is
naturally identified with the group
$\Gamma$~$\ref{The_group_Gamma}$, and we denote by
$$[x,u] \mapsto [x,u \#A] , \ \ \  A \in \Gamma$$ the action of
$\Gamma$ on $\widetilde {\cal L}$. Moreover, we denote by ${\cal
H}$ the set of all the zero-mean normalized Hamiltonian functions,
i.e.\
$${\cal H} =\Bigl \{ H \in C^{\infty}(S^1 \times M) \, \big | \, \int_M
H_t \, \omega^n =0, \ {\rm for \ all} \ t \in [0,1] \Bigr \}. $$
For $H \in {\cal H}$, the symplectic action functional ${\cal A}_H
\co \widetilde {\cal L} \to {\mathbb R}$ is defined as
$$ {\cal A}_H([x,u]) := - \int_u \omega + \int_{S^1} H(x(t),t)dt. $$
Note that $${\cal A}_H([x,u \# A]) = {\cal A}_H([x,u]) -
\omega(A).$$ Let ${\cal P}_H$ be the set of all contractible
$1$--periodic orbits of the Hamiltonian flow generated by $H$.
Denote by $\widetilde {{\cal P}_H }$ the subset of pairs $[x,u]
\in \widetilde {\cal L}$ where $ x \in {\cal P}_H$. It is not
difficult to verify that $\widetilde {{\cal P}_H}$ coincides with
the set of critical points of ${\cal A}_H$. We define the {\it
action spectrum} of $H$, denoted by ${\rm Spec}(H)$, as
$$ {\rm Spec}(H) := \left \{ {\cal A}_H(x,u) \in {\mathbb R} \ | \ [x,u] \in
\widetilde {{\cal P}_H} \right \}.$$ Recall that the action
spectrum is either a discrete or a countable dense subset of
${\mathbb R}$~\cite{Oh}.

We now turn to give the definition of the filtered Floer homology
group. For a generic $H \in {\cal H}$ and $\alpha \in  \{ {\mathbb
R} \setminus {\rm Spec}(H) \}  \cup \{ {\infty} \}$ define the vector
space $CF^{\alpha}_k(H)$ to be
$$  CF^{\alpha}_k(H) = \Bigl \{ \sum_{_{[x,u] \in \widetilde {\cal
P}(H)}} \beta_{[x,u]} [x,u]  \mid  \ \beta_{[x,u]} \in {\mathbb
Q}, \
  \mu([x,u]) = k, \  {\cal A}_H([x,u]) < \alpha
 \Bigr \},$$
where each sum satisfies the following finiteness condition:
$$ \# \left \{ [x,u] \in {\widetilde {{\cal P}_H}} \ | \ \beta_{[x,u]}
 \neq 0 \ {\rm and} \ {\cal A}_H([x,u]) > \delta \right \} < \infty,
\ {\rm for \ every} \ \delta \in {\mathbb R}.$$ Here $\mu([x,u])$
denotes the Conley--Zehnder index
 $\mu \co {\widetilde {{\cal P}_H}} \to {\mathbb Z}$
(see e.g.~\cite{Sa}) which satisfies $ \mu([x,u \# A]) -
\mu([x,u]) = 2c_1(A)$. In particular, the Conley--Zehnder index of
an element $x \in {\cal P}_H$ is well-defined modulo $2N_M$, where
$N_M$ is the minimal Chern number of $(M,\omega)$. The complex
$CF_k^{\infty}(H)$ is a module over the Novikov ring
$\Lambda$~$\ref{Def_Novikov_Ring}$, where the scalar
multiplication of $\xi \in CF_k^{\infty}(H)$ with $\lambda \in
\Lambda$ is given by
$$
 \sum_{A} \sum_{[x,u]} a_A \cdot \alpha_{[x,u]}  [x,u \# A]. $$
For each given $[x,w]$ and $[y,v]$ in $\widetilde {{\cal P}_H}$,
let ${\cal M}(H,J,[x,u],[y,v])$ be the moduli space of Floer
connecting orbits from $[x,w]$ to $[y,v]$, i.e.\ the set of
solutions $u \co {\mathbb R} \times S^1 \to M$ of the system
$$   \left \{ \begin{array}{c}
         \partial_{s}u + J_t(u)(\partial_t u - X_{H_t}(u)) = 0, \\
         \lim_{s \rightarrow - \infty} u(s,t) = x(t), \ \ \ \lim_{s \rightarrow \infty} u(s,t) = y(t), \\
       w \# u \#v \ {\rm represent \ the \ trivial \ class \ in \ } \Gamma.
 \end{array} \right \}   $$
It follows from the assumption of strongly semi-positivity and
from Gromov's compactness
 theorem~\cite{G} that for a generic choice of $J$ the moduli spaces ${\cal M}([x,u],
[y,v])$, for $\mu([x,u]) - \mu([y,v]) = 1$, are compact.

The Floer boundary operator $\partial \co CF_k^{\alpha}(H) \to
CF_{k-1}^{\alpha}(H)$ is defined by
$$ \partial([x,w]) = \sum n \bigl ([x,w],[y,v] \bigr )\ [y,v],$$
where the sum runs over all the elements $[y,v] \in {\widetilde
{{\cal P}_H}}$ such that $\mu[y,v] = k-1$ and $n \bigl
([x,w],[y,v] \bigr)$ denotes counting the (finitely many)
un-parameterized Floer trajectories with a sign determined by a
coherent orientation. As proved by Floer in~\cite{F}, the boundary
operator $\partial$ is well defined, satisfies $\partial^2 = 0$
and preserves the subspaces $CF_*^{\alpha}(H)$ (see~\cite{HS}).
Therefore, defining the quotient group by
$$ CF_*^{[a,b)}(H,J) = CF_*^b(H,J) \ / \ CF_*^a(H,J) \ \ \ \ \ \ \ \ \ (-\infty < a \leq b \leq \infty), $$
the boundary map induces a boundary operator $\partial \co
CF_*^{(a,b]}(H) \to CF_*^{(a,b]}(H)$, and we can define the Floer
homology group
 by $$ HF_*^{(a,b]}(J,H) =
(CF_*^{(a,b]}(H),\partial).$$
We will use the convention $HF_*(H,J) = HF_*^{ (-
\infty,\infty]}(H,J)$ and $HF_*^{a}(H,J) = HF_*^{(- \infty,a]}$.
The graded homology $HF_*(H,J)$ is a module over the Novikov ring
$\Lambda$, since the boundary operator is linear over $\Lambda$.
Note that these homology groups have been defined for generic
Hamiltonians only. However, one can extend the definition to all
$H \in {\cal H}$ using a continuation procedure (see
e.g~\cite{EP}). A key observation is that the Floer homology
groups are independent of the almost complex structure $J$ and the
Hamiltonian $H$ used to define them.
Moreover, if two Hamiltonian functions $H_1,H_2 \in {\cal H}$
generate the same element ${\varphi} \in \widetilde {\rm
Ham}(M,\omega)$, then ${\rm Spec}(H_1) = {\rm Spec}(H_2)$ (see~\cite{Oh2}
and~\cite{S}) and the spaces $ HF_*^{(a,b]}(J,H_1)$ and $
HF_*^{(a,b]}(J,H_2)$ can be canonically identified. Therefore, we
shall drop the notation $J$ and $H$ in $HF_*(H,J)$ and denote
$HF_*(\varphi) = HF_*(J,H)$ where $\varphi \in \widetilde {\rm
Ham}(M,\omega)$ is generated by $H$.

We denote by $\pi_{\alpha} \co HF_{*}({\varphi}) \to
HF^{(\alpha,\infty]}({ \varphi})$ the homomorphisms induced by the
natural projection $CF_{\infty}(H) \rightarrow
CF_{\infty}(H)/CF_{\alpha}(H)$ of Floer complexes and by $
i_{\alpha} \co HF_{\alpha}({ \varphi}) \to HF_{*}({ \varphi})$ the
homomorphism induced by the inclusion map $ i_{\alpha} \co
CF_*^{\alpha}(H) \to CF_*^{\infty}(H)$. Note that the homology
long exact sequence yields  ${\rm Kernel} \ \pi_{\alpha} = {\rm
Image} \ i_\alpha$. There exists a natural ring structure on the
Floer homology groups named {\it Pair-of-pants product} (see
e.g.~\cite{PSS}) $$ *_{{\rm pp}} \co HF_{\alpha}( \varphi) \times
HF_{\beta}( \psi) \to HF_{\alpha + \beta}({ \varphi \psi}).$$
In~\cite{PSS}, Piunikhin, Salamon and Schwarz constructed a
homomorphism between the Quantum homology groups $QH_*(M)$ and the
Floer homology groups $HF_*(M)$. Furthermore, they showed that the
homomorphism $ \Phi \co QH_*(M)
\to HF_*(H) $ 
is an isomorphism which preserves the grading and intertwines the
quantum product on $QH_*(M)$ with the pair-of-pants product on
$HF_*(H)$, i.e.\ $
 \Phi \left ( i_{\alpha + \beta } (\xi *_{pp} \eta) \right ) =
 \Phi \left (i_{\alpha}(\xi) \right ) * \Phi \left ( i_{\beta}(\eta) \right )$, for  every $\xi \in
 HF_{\alpha}(\varphi), \  \eta \in
  HF_{\beta}(\psi)$. In what follows, we will refer to the isomorphism $\Phi$ as the
PSS isomorphism.

\section{The existence of a Calabi quasi-morphism}

Let $(M^{2n},\omega)$ be a closed connected rational strongly
semi-positive symplectic manifold. Following the works of
Viterbo~\cite{V}, Schwarz~\cite{S}, and Oh~\cite{Oh1},
we recall the definition 
of a spectral invariant $c$ which plays a central role in the
proof of \fullref{The criterion}. We refer the readers
to~\cite{Oh1} and~\cite{MS2} for complete details of the
construction and proofs of the general properties of this spectral
invariant. A brief description of Floer homology and the PSS
isomorphism was also given in the above
\fullref{Pre_on_Floer}.

We define the spectral invariant $ c \co QH_*(M) \times
{\widetilde {\rm Ham}(M,\omega)} \to {\mathbb R} $ as follows. For
the elements $0 \neq a \in QH_*(M)$ and $\varphi \in {\widetilde
{\rm Ham}}(M,\omega)$, we set
$$ c(a, \varphi) = \inf \left \{ {\alpha} \in
{\mathbb R} \ | \ \Phi(a) \in {\rm Image} \  i_{\alpha} \right \},
$$ where $\Phi \co QH_*(M) \to HF_*(\varphi)$ is the PSS
isomorphism between the quantum homology and the Floer homology,
and $i_{\alpha} \co HF_{\alpha}(\varphi) \to HF_{*}(\varphi)$ is
the natural inclusion in the filtered Floer homology. The
non-trivial fact that $- \infty < c(a, \varphi) < \infty$ is
proved in~\cite{Oh1}. Moreover, $c(a,\varphi)$ has the following
properties~\cite{Oh1},~\cite{MS2}: For every $a,b \in QH_*(M)$ and
every $\varphi, \psi \in \widetilde {\rm Ham}(M)$

\begin{enumerate}
\item[(P1)]  $c(a*b,\varphi \psi) \leq c(a,\varphi) +
c(b,\psi)$, 
\item[(P2)]  $c(a,{\Id)} = val(a)$, 
\item[(P3)]  $c(a,\varphi) = \underset {m}
\sup \ c(a^{[m]},\varphi)$, 
\item[(P4)] $c(a  q^A,\varphi) = c(a,\varphi) + \omega(A)$, 
\ for every $q^A \in \Lambda$,
\end{enumerate}
where $a^{[m]}$ is the grade-$m$--component of $a$, $\Id$ is the
identity in ${\widetilde {\rm Ham}}(M,\omega)$ and $val(\cdot)$ is
the valuation function~$\ref{Def_valuation}$ defined in
\fullref{Sub-Sec_QH}.

The following lemma, which can be considered as a Poincar\'{e}
duality type lemma, enables us to compare the spectral invariants
of $\varphi$ and $\varphi^{-1}$. It is the analogue of Lemma 2.2
from \cite{EP} in the rational non-monotone case .

\begin{lemma} 
\label{PD_lemma} For every $0 \neq \gamma \in QH_*(M)$ and every
$\varphi \in { \widetilde {\rm Ham}} (M,\omega)$
$$ c(\gamma,\varphi) = - \inf \left \{  c(\delta,\varphi^{-1}) \ | \ \Pi(\delta,\gamma) \neq 0  \right
\}, $$ where $\Pi(\cdot,\cdot)$ is the ${\mathbb Q}$--valued
pairing~$\ref{Def_Poincare_pairing}$ defined in
\fullref{Sub-Sec_QH}.
\end{lemma}

The proof of the lemma is given in \fullref{Section_PD} below.
In order to prove \fullref{The criterion} we will also need
the following proposition. Assume that the subalgebra $QH_{2n}(M)
\subset QH_*(M)$ is semi-simple over the field $\Lambda_0$ and let
$QH_{2n}(M) = QH_{2n}^1(M) \oplus \cdots \oplus Q_{2n}^k(M)$ be a
decomposition of $QH_{2n}(M)$ into a direct sum of fields. Then we
have
\begin{proposition}  \label{bounded_eval_1}
There exists a positive constant $C \in {\mathbb R}$ such that for
every $0 \neq \gamma \in QH_{2n}^1(M)$
$$val(\gamma) + val(\gamma^{-1}) \leq C.$$
\end{proposition}

Postponing the proof of the above proposition we first present the
proof of \fullref{The criterion} and
\fullref{Main_Examples}. In the proof of \fullref{The
criterion} we follow the strategy of the proof used by Entov and
Polterovich in~\cite{EP}.

\begin{proof}[Proof of \fullref{The criterion}]
 Let $QH_{2n}(M) =
QH_{2n}^1(M) \oplus \cdots \oplus Q_{2n}^k(M)$ be a decomposition
of $QH_{2n}(M)$ into a direct sum of fields. Consider the map
$\widetilde{r} \co {\widetilde {\rm Ham}}(M) \to {\mathbb R}$
defined by:
$$ \ \ {\widetilde r}(\varphi) = -{\rm vol}(M) \cdot \lim_{n \rightarrow \infty}
{\frac {c(e_1,\varphi^n)}n},$$ where $e_1$ is the unit element of
$QH_{2n}^1(M)$. This is a standard homogenization of the map
$c(e_1,\cdot) \co {\widetilde {\rm Ham}}(M) \to {\mathbb R}$. We
claim that $\widetilde{r}$ is a Lipschitz homogenous Calabi
quasi-morphism. The proof of the Calabi property and the Lipschitz
property of $\widetilde r$ goes along the same lines as the proof
of Propositions $3.3$ and $3.5$ in~\cite{EP} with the notations
suitably adapted. Thus, we will omit the details of the proof of
these properties and concentrate on proving that $\widetilde r$ is
a quasi-morphism. We will show that $c(e_1,\cdot)$ is a
quasi-morphism, this immediately implies that its homogenization
$\widetilde r$ is also a quasi-morphism.

Notice that the upper bound follows easily from the triangle
inequality  (P$1)$:
$$ c \left (e_1,\varphi \psi \right ) = c \left (e_1 * e_1, \varphi  \psi \right)
\leq c \left (e_1,\varphi \right ) + c \left (e_1,\psi \right ).$$
Next, it follows from (P$1)$ and \fullref{PD_lemma} that
\begin{eqnarray*}  c \left (e_1,\varphi \right )  & \leq &
c \left ( e_1,\varphi  \psi \right ) + c  ( e_1,\psi^{-1}) =
 c \left (e_1,\varphi  \psi \right ) - \inf_{a:\Pi(a,e_1) \neq 0}
c \left ( a,\psi \right ).
\end{eqnarray*}
From the definition of the intersection pairing $\Pi$~$\ref{Def_Poincare_pairing}$ we have that 
$$\{a \ | \ \Pi(a,e_1) \neq 0 \} = \{a \ | \ \Pi(a^{[0]},e_1) \neq
0\} = \{a \ | \ \Pi(a^{[0]} * e_1,M) \neq 0 \}.$$ 
Combining this with the above property (P$3$) we may further
estimate
\begin{equation} \label{eq_1} c(e_1,\varphi)  \leq  c(e_1,\varphi  \psi) -
\inf_{a :  \Pi(a^{[0]} *e_1,M) \neq 0} c(a^{[0]},\psi).
\end{equation} Our next step is to find a lower bound for the term
$c(a^{[0]},\psi)$ provided that $ \Pi(a^{[0]} * e_1, M) \neq 0$.
For this, we shall first ``shift" and then ``project'', roughly
speaking, the element $a^{[0]} \in QH_0(M)$ to the field
$QH_{2n}^1(M)$. More precisely, since we assumed that the minimal
Chern number $N_M$ divides $n$, there exists an element $q^A$ in
the Novikov ring $\Lambda$ such that $a^{[0]} q^A \in QH_{2n}(M)$.
Thus, it follows from properties (P$1)$ and (P$4)$ that
\begin{equation} \label{equation5.2} c(a^{[0]},\psi) =
c(a^{[0]}q^A,\psi) - \omega(A) \geq c(e_1 * a^{[0]}q^A,\psi) -
c(e_1,{\mathbb \Id}) - \omega(A).
\end{equation} Moreover, it follows from the assumption 
$\Pi(a^{[0]} * e_1, M) \neq 0$ and from the definition of the
element $e_1$, that $ e_1 * a^{[0]}q^A \in QH_{2n}^1(M) \setminus
\{ 0 \}$. Hence, since $QH_{2n}^1(M)$ is a field, $ e_1 *
a^{[0]}q^A$ is an invertible
element inside it. 
Using the triangle inequality (P$1)$ once again we get
$$ c(e_1,\psi) \leq c(e_1 * a^{[0]}q^A,\psi) + c \bigl ( ( e_1 * a^{[0]}q^A)^{-1},{\mathbb \Id} \bigr ).$$
Here $( e_1 * a^{[0]}q^A)^{-1}$ is the inverse of $e_1 *
a^{[0]}q^A$ inside $QH_{2n}^1$. Next, by substituting this in the
above inequality~$\ref{equation5.2}$ and applying (P$2)$ we can
conclude
\[ c(a^{[0]},\psi)  \geq  c(e_1,\psi) -
 val \bigl ( ( e_1 * a^{[0]}q^A)^{-1} \bigr )  - val(e_1) - \omega(A). \]
By assigning this lower bound of $ c(a^{[0]},\psi)$
into~$\ref{eq_1}$ we further conclude
\begin{eqnarray*}
c(e_1,\varphi) & \leq  & c(e_1,\varphi \psi) -  c(e_1,\psi) +
\sup_{ {a: \Pi( a^{[0]} * e_1,M) \neq 0}} ~ val \left ( (e_1 *
a^{[0]} q^A)^{-1} \right ) + C',
\end{eqnarray*}
where $C'$ is the value $val(e_1) + \omega(A)$. The last step of
the proof is to find a universal upper bound for $val \bigl ( (e_1
* a^{[0]}q^A)^{-1} \bigr )$ provided that $\Pi(a^{[0]} * e_1,M)
\neq 0$. Note that the condition $\Pi(a^{[0]} * e_1,M) \neq 0$
 implies that $val \left (e_1 * a^{[0]} \right ) \geq 0$ and hence $val
\left (e_1 * a^{[0]}q^A \right ) \geq -\omega(A)$. Therefore, 
it follows from \fullref{bounded_eval_1} 
that
$val \bigl ( (e_1
* a^{[0]} q^A)^{-1} \bigr ) <
C + \omega(A)$. We have shown that $c(e_1,\cdot)$ is a
quasi-morphism, the proof of the theorem is thus complete.
\end{proof}

\begin{proof}[Proof of \fullref{Main_Examples}]
Let $M$ be one of the manifolds $(X_{\lambda},\omega_{\lambda})$
or $(Y_{\mu},\omega_{\mu})$ listed in the theorem. 
It follows from \fullref{QH_4_is_a_subfield_1} and
\fullref{QH_4_is_a_subfield_2} that the subalgebra $QH_{4}(M)$
is semi-simple. Moreover, the minimal Chern number of
$(X_{\lambda},\omega_{\lambda})$ and $(Y_{\mu},\omega_{\mu})$ is
$2$ and $1$ respectively. Thus, it follows from \fullref{The
criterion} that there exists a Lipschitz homogeneous Calabi
quasi-morphism $r \co {{\widetilde {\rm Ham}}(M,\omega)} \to
{\mathbb R}$ as required.
\end{proof}

\begin{remark} {\rm
As mentioned in \fullref{algebraic_structure_of_QH(Y)} above,
in the case of $(Y_{\mu},\omega_{\mu})$ where $ 0 < \mu < {\frac 1
3} $, the subalgebra $QH_4(Y_{\mu})$ splits into direct sum of two
fields. Thus, using the units of these fields alternately,
\fullref{The criterion} implies the existence of two Calabi
quasi-morphisms. We do not know whether they are equivalent or
not.}
\end{remark}

We return now to the proof of \fullref{bounded_eval_1}. We
will follow closely Lemma 3.2 in~\cite{EP}.
\begin{proof}[Proof of \fullref{bounded_eval_1}]
From the definition of the graded Novikov ring it follows that
$\Lambda_{0}$ can be identified with the field ${\cal R} =
{\mathbb Q}[[x]$ of Laurent series $\sum \alpha_j x^j$ with
coefficients in ${\mathbb Q}$ and $\alpha_j =0$ for large enough
$j$'s. Moreover, it is not hard to check that $QH_{k}(M)$ is a
finite dimensional module over $\Lambda_{0}$. We denote by $\sigma
\co {\cal R} \to {\mathbb Z}$ the map which associates to a
nonzero element $ \sum \alpha_j x^j \in {\cal R}$ the maximal $j$,
such that $\alpha_j \neq 0$. We set $\sigma(0) = - \infty$. For
$\kappa \in {\cal R}$, put $| \kappa |_1 = {\rm exp} \
\sigma(\kappa)$. Thus, $| \cdot |_1$ is a non-Archimedean absolute
value on ${\cal R}$ and moreover, ${\cal R}$ is complete with
respect to $| \cdot |_1$.  For preliminaries on non-Archimedean
geometry we refer the readers to~\cite{Ar}. Since the field
$QH_{2n}^1(M)$ can be considered as a finite dimensional vector
space over ${\cal R}$, the absolute value $| \cdot |_1$ can be
extended to an absolute value $| \cdot |_2$ on $QH_{2n}^1(M)$
(see~\cite{Ar}). Note that $| \cdot |_2$ induces a multiplicative
norm $\| \cdot \|_2$ on $QH_{2n}^1(M)$. On the other hand, we can
consider a different norm on $QH_{2n}^1(M)$ defined by $\| \gamma
\|_3 = {\rm exp} \ val(\gamma)$. Since all the norms on a finite
dimensional vector space are equivalent, there is a constant $ C_1
> 0$ such that
$$\| \gamma \|_3 \leq C_1 \cdot \| \gamma \|_2, \ {\rm for \ every \ } 0 \neq \gamma \in QH_{2n}^1(M). $$
Hence, for $ 0 \neq \gamma \in QH_{2n}^1(M)$, we have
$$ \| \gamma \|_3 \cdot \| \gamma^{-1} \|_3 \leq  C_1^2 \cdot \| \gamma \|_2 \cdot \| \gamma^{-1} \|_2 = C_1^2. $$
Therefore, $val(\gamma) + val(\gamma^{-1}) \leq C$ where $ C = 2
\log C_1$. This completes the proof of the proposition.
\end{proof}

\section[Restricting  r to the fundamental group of Ham(M)]{Restricting $\widetilde r$  to the fundamental group of $ {\rm Ham}(M)$ } \label{The_restriction_to_Ham}
In this section we discuss the restriction of the above mentioned
Calabi quasi-morphism $\widetilde r \co {{\widetilde {\rm
Ham}}(M)} \to {\mathbb R}$, where $M$ is one of the manifolds
listed in \fullref{Main_Examples}, to the abelian subgroup
$\pi_{1} \bigl ( {\rm Ham(M)} \bigr) \subset {{\widetilde {\rm
Ham}}(M)}$. 
For this purpose, we follow~\cite{EP} and use the
 Seidel representation $\Psi \co
\pi_1 \bigl ( {\rm Ham}(M) \bigr ) \to QH_{ev}^{\times}
(M,{\mathbb R})$ (see e.g.~\cite{Se}, \cite{LMP}), where
$QH_{ev}^{\times}(M,{\mathbb R})$ denotes the group of units in the
even part of the quantum homology algebra of $M$ with coefficients in
a real Novikov ring.  We start with the following preparation.

\subsection{Hamiltonian fibrations over the two sphere }
There is a one-to-one correspondence between homotopy classes of
loops in Ham$(M)$ and isomorphism classes of Hamiltonian
fibrations over
the two-sphere $S^2$ 
given by the following ``clutching" construction (see
e.g.~\cite{Se},~\cite{LMP}). We assign to each loop $\varphi = \{
\varphi_t \} \in {\rm Ham}(M)$ the bundle $(M,\omega) \rightarrow
P_{\varphi} \rightarrow S^2$ obtained by gluing together the
trivial fiber bundles $D^{\pm} \times M$ along their boundary via
$(t,x) \mapsto (t,\varphi_t(x))$. Here we consider $S^2$ as $D^{+}
\cup D^{-}$, where $D^{\pm}$ are closed discs with boundaries
identified with $S^1$. Moreover, we orient the equator $D^{+} \cap
D^{-}$ as the boundary of $D^{+}$. Note that this correspondence
can be reversed.

As noted in~\cite{Se}, there are two canonical cohomology classes
associated with such a fibration. One is the {\it coupling class}
$u_{\varphi} \in H^2(P_{\varphi},{\mathbb R})$ which is uniquely
defined by the following two conditions: the first is that it
coincides with the class of the symplectic form on each fiber, and
the second is that its top power $u_{\varphi}^{n+1}$ vanishes. The
other cohomology class is the first Chern class of the vertical
tangent bundle $c_{\varphi} = c_1(TP^{\rm vert}_{\varphi}) \in
H^2(P_{\varphi},{\mathbb R}).$ We define an equivalent relation on
sections of the fibration $P_{\varphi} \rightarrow S^2$ in the
following way: First, equip $S^2$ with a positive oriented complex
structure $j$, and $P_{\varphi}$ with an almost complex structure
$J$ such that the restriction of $J$ on each fiber is compatible
with the symplectic form on it, and the projection $\pi \co
P_{\varphi} \to S^2$
is a $(J,j)$--holomorphic map. 
Next, two $(J,j)$--sections, $\nu_1$ and $\nu_2$, of $\pi \co
P_{\varphi} \to S^2$ are said to be $\Gamma$--{\it equivalent} if
$$u_{\varphi}[\nu_1(S^2)] = u_{\varphi}[\nu_2(S^2)], \ \
c_{\varphi}[\nu_1(S^2)] = c_{\varphi}[\nu_2(S^2)].$$ It has been
shown in~\cite{Se} that the set ${\cal S}_{\varphi}$ of all such
equivalent classes is an affine space modeled on the
 group $\Gamma$~$\ref{The_group_Gamma}$.

\subsection{The Seidel representation}

The following description of the Seidel representation, which is
somehow different from Seidel's original work, can be found
in~\cite{LMP}.
For technical reasons, 
 it will be more convenient to work in what follows with a
slightly larger Novikov ring than in~$\ref{Def_Novikov_Ring}$.
More precisely, set ${\cal H}_{\mathbb R} := H_2^S(M,{\mathbb R})
~ / ~ ({\rm ker} \, c_1 \, \cap \, {\rm ker} \, \omega) $, where
$H_2^S(M,{\mathbb R})$ is the image of $\pi_2(M)$ in
$H_2(M,{\mathbb R})$. We define the {\it real Novikov ring} as
$$
 \Lambda_{\mathbb R} = \Bigl \{ \sum_{A \in {\cal H}_{\mathbb R}} \lambda_A q^A \ | \
 \lambda_A \in {\mathbb Q}, \  \#  \{ A \in \Gamma \ | \ \lambda_A \neq 0, \ \omega(A) >  c \}
 < \infty, \ \forall c \in {\mathbb R}\, \Bigr \}, $$
and set 
$QH_*(M):= QH_*(M,\Lambda_{\mathbb R}) = H_*(M) \otimes
\Lambda_{\mathbb R}$ to be the real quantum homology of $M$.

Next, let $\varphi$ be a loop of Hamiltonian diffeomorphisms and
$\nu$ be an equivalence class of sections of $P_{\varphi}$. Set
$d=2c_{\varphi}(\nu)$. We define a $\Lambda_{\mathbb R}$--linear
map $\Psi_{\varphi,\nu} \co QH_*(M) \to QH_{* + d}(M)$ as follows:
for $a \in H_*(M,{\mathbb Z})$, $\Psi_{\varphi,\nu}$ is the class
in $QH_{*+d}(M)$ whose intersection with $b \in H_*(M,{\mathbb
Z})$ is given by
$$ \Psi_{\varphi,\nu}(a) \cdot_M b = \sum_{B \in {\cal H}}n_{P_{\varphi}}(i(a),i(b);\nu + i(B)) q^{-B},$$
where $i$ is the homomorphism $H_*(M) \rightarrow
H_*(P_{\varphi})$, the intersection $\cdot_M$ is the linear
extension to $QH_*(M)$ of the standard intersection pairing on
$H_*(M,{\mathbb Q})$, and $n_{P_{\varphi}}(v,w;\mu)$ is the
Gromov--Witten invariant which counts isolated $J$--holomorphic
stable curves in $P_{\varphi}$ of genus $0$ and two marked points,
such that each curve represents  the equivalence class $\mu$ and
whose marked points go through given generic representatives of
the classes $v$ and $w$ in $H_*(P_{\varphi},{\mathbb Z})$. When
the manifold $M$ is strongly semi-positive, these invariants are
well defined. Moreover, it follows from Gromov's compactness
theorem (see~\cite{G}) that for each given energy level $k$, there
are only finitely many section-classes $\mu = \nu + i(B)$ with
$\omega(B) \leq k$ that are represented by a $J$--holomorphic curve
in $P_{\varphi}$. Thus, $ \Psi_{\varphi,\nu}$ satisfies the
finiteness condition for elements in $QH_*(M)$.

For reasons of dimension, $n_{P_{\varphi}}(v,w;\mu)= 0 $ unless
$2c_{\varphi}(\mu) + {\rm dim}(v) + {\rm dim}(w) = 2n$. Thus,
$$ \Psi_{\varphi,\nu}(a) = \sum  a_{\nu,B} \, q^{-B}, \ \ a_{\nu,B}
\in H_*(M),$$ where $a_{\nu,B} \cdot_M b =
n_{P_{\varphi}}(i(a),i(b);\nu+i(B))$, and
$$ {\rm dim}(a_{\nu,B}) = {\rm dim}(a) + 2c_{\varphi}(\nu + i(B)) = {\rm dim}(a) +
2c_{\varphi}(\nu)+ 2c_{1}(B).$$ Note also that $ \Psi_{\varphi,\nu
+ A} = \Psi_{\varphi,\nu} \otimes q^A$. It has been shown by
Seidel~\cite{Se} (see also~\cite{LMP}) that $\Psi_{\varphi,\nu}$
is an isomorphism for all loops $\varphi$ and sections $\nu$.

Next, we use $\Psi_{\varphi,\nu}$ to define the Seidel
representation $$\Psi \co \pi_1({\rm Ham}(M)) \to
QH_*(M,\Lambda_{\mathbb R})^{\times}.$$ In order to do so, we take
a canonical section class $\nu_{\varphi}$ that (up to equivalence)
satisfies the composition rule $\nu_{\varphi  \psi} = \nu_{\psi}
\# \nu_{\varphi}$ , where $\nu_{\varphi  \psi}$ denotes the
obvious union of the sections in the fiber sum $P_{\psi  \varphi}
= P_{\psi} \# P_{\varphi}$. The section $\nu_{\varphi}$ is
uniquely determined by the requirement that \begin{equation*}
\label {normalization-of-the-section} u_{\varphi}(\nu_{\varphi})=
0 \ \ {\rm and} \ \ c_{\varphi}(\nu_{\varphi}) = 0. \end{equation*}
Moreover, it satisfies the above mentioned composition rule.
Therefore, we get a group homomorphism
$$ \rho \co \pi_1({\rm Ham}(M)) \to {\rm
Hom}_{\Lambda_{\mathbb R}}(QH_*(M,\Lambda_{\mathbb R})).$$ It has
been shown in~\cite{Se} that for all $\varphi \in \pi_1({\rm
Ham}(M))$ we have $$ \rho(\varphi)(a) =
\Psi_{\varphi},\nu_{\varphi}([M]) *_M a .$$ The Seidel
representation is defined to be 
the natural homomorphism
$$ \Psi \co \pi_1( {\rm Ham}(M) ) \to QH_*(M,
\Lambda_{\mathbb R})^{\times},$$ given by $\varphi \mapsto
\rho(\varphi)([M])$.

\subsection{Relation with the spectral invariant}
Throughout, $\pi_1({\rm Ham}(M))$ is considered as the group of
all loops in ${\rm Ham}(M)$ based at the identity $\Id \in {\rm
Ham}(M)$. Let ${\cal L}$ be the space of all smooth contractible
loops $x \co S^1 = {\mathbb R} / {\mathbb Z}
\to M$ and $\widetilde {\cal L}$ its cover 
introduced in Section $4$. Let $\varphi$ be a loop of Hamiltonian
diffeomorphisms.  It is known (see e.g.~\cite{Se}) that the orbits
$\varphi_t(x)$ of $\varphi$ are contractible. We consider the map
$T_{\varphi} \co {\cal L} \to {\cal L}$ which takes the loop
$x(t)$ to $\varphi_t(x(t)).$ 
In~\cite{Se}, Seidel showed that this action can be lifted (not
uniquely) to $\widetilde {\cal L}$. In fact it is not hard to
check that there is a one-to-one correspondence between such lifts
of $T_{\varphi}$ and equivalence classes of sections $\nu \in
{\cal S}_{\varphi}$. We denote by $\widetilde T_{\varphi,\nu}$ the
lift corresponding to $\nu \in {\cal S}_{\varphi}$. Next, let
$\varphi \in \pi_1({\rm Ham}(M))$ be a given loop generated by a
normalized
Hamiltonian $K 
\in {\cal H}$. 
The following formula, which can be found in~\cite{Se}
and~\cite{LMP},
 enables us to relate the Seidel representation with the spectral
invariant $c$ used to define the Calabi quasi-morphism $\widetilde
r$:
\begin{equation} \label{equation_about_action_shift}
 (\widetilde T_{\varphi,\nu}^*)^{-1} {\cal A}_H - {\cal A}_{K
 \sharp H} = - u_{\varphi}(\nu), \ \ {\rm for \ every} \ H \in
 {\cal H}.
\end{equation}
It has been shown in~\cite{Se} that the isomorphism in the quantum
homology level described in Section $6.2$, which is
obtained by multiplication with
$\Psi_{\varphi},\nu_{\varphi}([M])$ corresponds, under the
identification between the Floer and the quantum homology, to the
isomorphism $i \co HF_{\alpha}(H) \to HF_{\alpha +
u_{\varphi}(\nu) }(K \# H)$ induced by the action of $(\widetilde
T_{\varphi,\nu}^*)$ on $\widetilde {\cal L}$. The following
proposition can be found in~\cite{Oh3} or~\cite{EP}.

\begin{proposition} \label{formula-for-restriction}
For every loop $ [ \varphi ] \in \pi_1({\rm Ham}(M)) \subset
{\widetilde {\rm Ham}(M)} $ and every $a \in QH_*(M)$ we have $$
c(a,[\varphi]) = val(a * \Psi([{\varphi}])^{-1}). $$
\end{proposition}

\begin{proof}
Let $K \in {\cal H}$ be the normalized Hamiltonian function
generating the loop $[\varphi]$, and let $H \in {\cal H}$ be the
zero Hamiltonian generating the identity.
 The proposition immediately follows
from~$\ref{equation_about_action_shift}$ applied to $H$ and $K$.
\end{proof}

\section[Proof of \ref{Theorem-Restriction}]{Proof of \fullref{Theorem-Restriction}}

Recall that a homogeneous quasi-morphism on an abelian group is
always a homomorphism (see e.g.~\cite{EP}). Hence, in order to
prove \fullref{Theorem-Restriction}, we need to show that for
the manifolds listed in the theorem, the restriction of the Calabi
quasi-morphism $\widetilde r \co {{\widetilde {\rm Ham}}(M)} \to
{\mathbb R}$ on the fundamental group of ${\rm Ham}(M)$ is
non-trivial. We will divide the proof into two parts.

\subsection{ The case of $S^2 \times S^2$.}

Let $X_{\lambda} = S^2 \times S^2$ be equipped with the split
symplectic form $\omega_{\lambda} = \omega \oplus \lambda \,
\omega$, where $1 < \lambda $. As mentioned in Section $1$, 
there is an element $[\varphi]$ of infinite order in the
fundamental group of ${\rm Ham}(X_{\lambda})$ (see~\cite{Mc}).
This element can be represented by the following loop of
diffeomorphisms
$$ \varphi_t(z,w) = \bigl (z, \Upsilon_{z,t}(w) \bigr ),$$
where $\Upsilon_{z,t}$ denotes the
 $2 \pi t$--rotation of
the unit sphere $S^2$ around the axis through the points $z, \
-z$. Seidel showed in~\cite{Se} (see also~\cite{McTol}), by direct
calculation, that
$$\Psi([\varphi])^{-1}=(A-B)\otimes  q^{\alpha A + \beta B} \bigl ( \sum_{j=0}^{\infty} q^{j(A-B)} \bigr ),
$$
where $A$ and $B$ in $H_*(X_{\lambda})$  are the classes of  $[S^2
\times {\rm point}]$ and $[{\rm point} \times S^2]$ respectively,
and $\alpha,\beta \in {\mathbb R}$ were chosen such that $2c_1
(\alpha A + \beta B)=1$ and $\omega_{\lambda}(\alpha A + \beta B)=
{\frac 1 2} + {\frac 1 {6{\lambda}}}$.

\begin{lemma} \label{val_first_case} For every $n \in {\mathbb N}$ we have
$$val(\Psi([\varphi])^{-2n}) = 1 +  {\frac n {3\lambda}}.$$
\end{lemma}
\begin{proof}
First note that $val \bigl ((A-B)^{2n} \bigr ) = \max \{
val(A^kB^{2n-k}) \}$, where $0 \leq k \leq 2n$. Next, set
$\alpha_k = A^kB^{2n-k}$. It follows from the quantum
multiplication
relations~$\ref{quantum-relations-of-the-first-example}$ that
$val(\alpha_{k + 2}) = val(\alpha_k) + (\lambda -1) $ for every $0
\leq k \leq 2n-2$. Thus, $$val \bigl ((A-B)^{2n} \bigr ) = \max \{
val(A^{2n}) , val(A^{2n-1}B)  \} = -n+1.$$ Set $\Delta= q^{\alpha
A + \beta B} \bigl ( \sum_{j=0}^{\infty} q^{j(A-B)} \bigr )$. It
follow immediately that
$$ val(\Delta^{2n}) = val \bigl( q^ {2n (\alpha A + \beta B)}  \bigr) =
2n \Bigl  ( {\frac 1 2} + {\frac 1 {6{\lambda}}} \Bigr )= n +
{\frac n {3\lambda}}.
$$
This completes the proof of the Lemma.
\end{proof}

It follows from \fullref{QH_4_is_a_subfield_1} that the
subalgebra $QH_4(X_{\lambda})$ is a field. Thus, combining
\fullref{formula-for-restriction} and
\fullref{val_first_case}, we conclude that$$ {\widetilde
r}([\varphi]) =
  -{\rm vol}(X_{\lambda}) \cdot \lim_{n \rightarrow \infty} {\frac {val
\left (\Psi([{\varphi}])^{-2n} \right )} {2n}} =   - {\frac
{1+\lambda} {6\lambda}} \neq 0.
$$
We have shown that the restriction of the quasi-morphism
$\widetilde r \co {{\widetilde {\rm Ham}}(X_{\lambda})} \to
{\mathbb R}$
on the fundametal group of ${\rm Ham}(X_{\lambda})$ is non-trivial.
 This concludes the proof of
\fullref{Theorem-Restriction} for the above case.

\subsection{ The case of ${\mathbb C}P^2 \# \overline{{\mathbb
C}P^2}$}

Let $ Y_{\mu}=  {\mathbb C}P^2 \# {\overline {{\mathbb C}P^2}}$ be
the symplectic one-point blow-up of ${\mathbb C}P^2$ introduced in
Section $1$, equipped with the symplectic form $\omega_{\mu}$,
where ${\frac 1 3} \neq \mu \in (0,1)$. We will use here the same
notation as in Subsection $3.3$. It has been shown by Abreu--McDuff
in~\cite{AM} that the fundamental group of ${\rm Ham}(Y_{\mu})$ is
isomorphic to ${\mathbb Z}$ with a generator given by the rotation
$$ \varphi \, \co \, (z_1,z_2) \to (e^{-2\pi i t}z_1,z_2), \ \ \ 0 \leq
t \leq 1.$$
The Seidel representation of $\varphi$ was computed
in~\cite{Mc1},~\cite{McTol}  to be \begin{equation}
\label{Seidel-element-blow-up} \Psi([\varphi])^{-1} = P \otimes
q^{E/2 + 3F/4 - \delta(F-2E)}, \ \ \ {\rm where} \ \ \delta =
{\frac {(1-\mu)^2} {12 (1+ \mu) (1-3\mu) } }. \end{equation}

The following lemma can be immediately deduced from Lemma $5.1$
and Remark $5.5$ which both appear in~\cite{Mc1}.

\begin{lemma}  Let
${\frac 1 3} \neq \mu \in (0,1)$. Then
\begin{equation} \label{Val-of-Seidel-element-lemma} \lim_{k \rightarrow \infty} {\frac {
val(\Psi([\varphi])^{-k}) } k} = \left \{
\begin{array}{ll}
-\delta \, \omega(F-2E) ,  \ \  & {\frac 1 3} < \mu < 1 \\ \\
{\frac {12 - \delta} {12}} \,  \omega(F-2E), \ \ & 0 < \mu <
{\frac 1 3}.
\end{array} \right. \end{equation}
\end{lemma}

\begin{proof}
Denote by $Q$ the element $ P \otimes q^{E/2+3F/4}$ and consider
its powers $Q^k$ where $k \in {\mathbb N}$. It follows from the
quantum multiplicative relations discussed in Subsection $3.2$
that the only two possible cycles obtained by multiplication by
$Q$ are
$$ P \otimes q^{E/2+3F/4} \rightarrow E \otimes q^{F/2}
\rightarrow F\otimes q^{E/2+F/4} \rightarrow M \rightarrow P
\otimes q^{E/2+3F/4},$$ and
$$ P \otimes q^{E/2+3F/4} \rightarrow
F \otimes q^{F/2} \rightarrow M \otimes q^{F/4-E/2} \rightarrow P
\otimes q^{F}.$$
Thus, since the first cycle does not change the valuation, while
the second cycle increases it 
by $\omega(F/4 - E/2)$, we have that $val(Q^k)$ is either bounded
as $k \rightarrow \infty$ when $\omega(F/4-E/2) < 0$ or linearly
grows otherwise. Hence, we get that
$$
val(\Psi([\varphi])^{-k})   = \left \{ \begin{array}{ll}
C -\delta \, k \, \omega(F-2E) ,  \ \  & {\frac 1 3} < \mu < 1 \\ \\
C + {\frac k 3 } \, {\omega(F/4-E/2)} - \omega(\delta(F-2E)),  \ \
& 0 < \mu < {\frac 1 3},
\end{array} \right. $$
where $C$ is some universal constant. This completes the proof.
\end{proof}

A straightforward calculation shows that the above
expression~$\ref{Val-of-Seidel-element-lemma}$ is strictly
negative for every $0 < \mu < 1$. Thus, it follows from
\fullref{formula-for-restriction} and the fact that
$val(a*b) \leq val(a) + val(b)$ that
$$ {\widetilde r}([\varphi]) =
  -{\rm vol}(Y_{\mu})  \lim_{k \rightarrow \infty} {\frac {val
\left (e_1 * \Psi([{\varphi}])^{-k} \right )} {k}}  \geq  -{\rm
vol}(Y_{\mu})  \lim_{k \rightarrow \infty} {\frac {val
\left ( \Psi([{\varphi}])^{-k} \right )} {k}}   > 0 $$
Hence, the restriction of the quasi-morphism $\widetilde r \co
{{\widetilde {\rm Ham}}(Y_{\mu})} \to {\mathbb R}$ on the
fundamental group of ${\rm Ham}(Y_{\mu})$ is non-trivial. The
proof of \fullref{Theorem-Restriction} is now complete.

\section[Proof of \ref{PD_lemma}]{Proof of \fullref{PD_lemma}} \label{Section_PD}

Let $(M,\omega)$ be a closed, rational and strongly semi-positive
symplectic manifold of dimension $2n$. 
Note that for rational symplectic manifolds the action spectrum is
a discrete subset of ${\mathbb R}$, and thus there are only a
finite number of critical values of the action functional ${\cal
A}_{H}$ in any finite segment $[a,b] \subset {\mathbb R}$. Let
$H(t,x) \in {\cal H}$ be a Hamiltonian function generating
$\varphi \in \widetilde {{\rm Ham}}(M)$, and denote by $
{\widetilde H}(t,x) = -H(-t,x)$ the Hamiltonian function
generating the inverse symplectomorphism $\varphi^{-1}$. The set
${\cal P}_{H}$ of critical points of ${\cal A}_H$ is isomorphic to
${\cal P}_{\widetilde {H}}$ via ${\widetilde x}(t) = x(-t)$, and
$[x,u] \in {\widetilde {{\cal P}_H}}$ corresponds to
$[\widetilde{x},\widetilde{u}] \in {\widetilde {{\cal
P}_{\widetilde {H}}}}$ where
$$ {\widetilde u}(s,t)
= u(-s,-t), \  \ \mu([{\widetilde x},{\widetilde u}]) = 2n -
\mu([x,u]) \ \ {\rm and} \ \  {\cal A}_H([x,u]) = - {\cal
A}_{\bar{H}}([\widetilde {x},\widetilde {u}]). $$
We define a pairing $L \co CF_{k} (H) \times CF_{2n-k}(\widetilde
H) \to \Lambda_0$ by
\begin{equation} \label{The-Pairing-L}  L \left ( \sum
\alpha_{[x,u]} \cdot [x,u] , \sum \beta_{[\widetilde y,\widetilde
v]} \cdot [\widetilde y,\widetilde v] \right  ) = \sum_{A} \Bigl (
\sum_{[x,u]} \alpha_{[x,u]} \cdot \beta_{[x,u \sharp -A]} \Bigr )
q^A,
\end{equation}
 where the inner sum runs over all pairs $ [x,u] \in
{\widetilde {{\cal P}_{H}}}$ and the outer sum runs over all $A
\in \Gamma$ with $c_1(A) = 0$. The pairing $L$ is well defined.
Indeed, consider first the inner sum, the finiteness condition in
the definition of $CF_*(H)$ implies that it contains only finitely
many elements. Secondly, it follows from the same reason that the
power series on the right hand side of~$\ref{The-Pairing-L}$
satisfies the finiteness condition from the definition of the
Novikov ring~$\ref{Def_Novikov_Ring}$. It is not hard to check
that the pairing $L$ is linear over $\Lambda_0$ and that it is
non-degenerate in the standard sense.
Thus, since the vector spaces $CF_{k}(H)$ and
$CF_{2n-k}(\widetilde H)$ are finite dimensional over $\Lambda_0$,
which is in our case a field, the pairing $L$ determines an
isomorphism
$$ CF_k(H) \simeq {\rm Hom}_{\Lambda_0} \left ( CF_{2n-k}(\widetilde
H),\Lambda_0 \right ). $$
From the universal coefficient theorem we obtain a Poincar\'{e}
duality isomorphism $$  HF_k(H)   \simeq {\rm Hom}_{\Lambda_0}
\left ( HF_{2n-k}(\widetilde H),\Lambda_0 \right ) .$$
In~\cite{PSS} it has been shown that the pairing determined by
this isomorphism, which by abuse of notation we also denote by
$L$, agrees with the intersection pairing $\Delta(\cdot,\cdot)$ on
the quantum homology $QH_*(M)$. More precisely, let $\Phi \co
QH_*(M) \to HF_*(H)$ be the PSS isomorphism described in Section
$4$. Then, for every $a \in HF_k(H)$ and $b \in QH_{2n-k}(M)$ we
have \begin{equation} \label{CDF1} \Delta  ( \Phi^{-1}(a),b ) = L
\left ( a, \Phi(b) \right ). \end{equation}
Next, we consider the filtered Floer homology complexes $CF_k^{(-
\infty, \alpha]}(H)$ and\break $CF_{2n-k}^{(- \alpha,\infty]}(\widetilde
H)$. Note that these spaces are no longer vector spaces over
$\Lambda_0$ since they are not closed with respect to the
operation of multiplication by a scalar. We define a ${\mathbb
Q}$--valued pairing $L' \co CF_k^{[- \infty,\alpha]}(H) \times
CF_{2n-k}^{(\alpha, \infty]}(\widetilde H) \to {\mathbb Q}$ by
$$ L' \left ( \sum \alpha_{[x,u]} \cdot [x,u],\sum \beta_{[\widetilde
y,\widetilde v]} \cdot [\widetilde y,\widetilde v] \right ) =
\sum_{[x,u]}(\alpha_{[x,u]} \cdot \beta_{[x,u]}).$$ This pairing
is well defined since any element in $CF_{2n-k}^{(-
\alpha,\infty]}(\widetilde H)$ is a finite sum. It is
straightforward to check that the pairing $L'$ is non-degenerate
in the standard sense 
and that it coincides with the zero term of $L$. In other words,
denote by $\tau \co \Lambda_0 \to {\mathbb Q}$ the map sending
$\sum a_A q^A$ to $a_0$, then for any $a \in CF_k^{(-
\infty,\alpha]}(H)$ and $b \in CF_{2n-k}(\widetilde H)$ we have
\begin{equation} \label{CDF2} \tau L \left (i_{\alpha,H}(a) , b
\right ) = L' (a,\pi_{- \alpha,\widetilde H}(b) ).\end{equation}
By abuse of notation, we also denote by $L'$ the induced pairing
in the homology level: $L' \co HF_k^{(- \infty,\alpha]}(H) \times
HF_{2n-k}^{(- \alpha,\infty]}(\widetilde H) \to {\mathbb Q}$.
Next, consider the following diagram:
$$ \begin{CD}
{{QH_k(M)}} @<{\Phi^{-1}}<< HF_k(H) @<{i_{{\alpha},H}}<<
 {HF_k^{(-\infty,\alpha]}( H)} \\
\times_{\Pi} @. \times_{\tau L} @. \times_{L'} \\
QH_{2n-k}(M) @>{\Phi}>> HF_{2n-k}(\widetilde H) @>{\pi_{{-
\alpha},\widetilde H} }>>
HF_{2n-k}^{(-{\alpha},\infty]}(\widetilde H) \\
@VVV @VVV  @VVV \\
{\mathbb Q} @. {\mathbb Q} @. {\mathbb Q}
\end{CD}
$$
\pagebreak

Combining equations~$\ref{CDF1}$ and~$\ref{CDF2}$ together we
conclude that
for every element $a \in HF_k^{(- \infty,\alpha]}(H)$ and $b \in
QH_{2n-k}(M)$ we have
\begin{equation} \label{CDF3} \Pi ( \Phi^{-1} \circ i_{{\alpha},H}(a),b
 ) = \tau L ( i_{{\alpha},H}(a),\Phi(b))
 = L' ( a,\pi_{{- \alpha},\widetilde H} \circ
\Phi(b) )
\end{equation}
We are now in a position to prove \fullref{PD_lemma}.
\begin{proof}[Proof of \fullref{PD_lemma}]
We divide the proof into two steps.
\begin{enumerate}
\item 
Fix an arbitrary $\varepsilon > 0$ and put $\alpha = \varepsilon -
c(\gamma,\varphi)$. It follows from the definition of the spectral
invariant $c$ 
 that $\Phi (
\gamma) \notin {\rm Image} \ i_{- \alpha, \varphi} $. Note that
${\rm Image} \ i_{- \alpha,\varphi} = {\rm Kernel} \ \pi_{-
\alpha,\varphi}$ and thus 
$\xi := \pi_{- \alpha,\varphi} \circ \Phi (\gamma) \neq 0$. Since
the pairing $L'$ is non-degenerate there exists $\eta \in
HF_{2n-*}^{(-\infty,\alpha]}(\varphi^{-1})$ such that
$L'(\eta,\xi) \neq 0$. From~$\ref{CDF3}$ we have that
$\Pi(\delta_0,\gamma) \neq 0$, where $\delta_0 = \Phi^{-1} \circ
i_{\alpha,\varphi^{-1}}(\eta)$. It follows from the definition
that $c(\delta_0,\varphi^{-1}) \leq \alpha$ and hence
$$ \inf_{_{\delta : \Pi(\delta,\gamma) \neq 0}} c(\delta,\varphi^{-1})
\ \leq c(\delta_0,\varphi^{-1}) \leq \alpha = \varepsilon -
c(\gamma,\varphi)$$ This inequality holds for every $\varepsilon >
0$, hence we conclude that $$\inf_{_{\delta : \Pi(\delta,\gamma)
\neq 0}} c(\delta,\varphi^{-1}) \leq - c(\gamma,\varphi).$$ \item
Fix an arbitrary $\varepsilon > 0$ and put $\alpha = - \varepsilon
- c(\gamma,\varphi)$. From the definition of $c(\cdot,\cdot)$ it
follows that $\Phi ( \gamma) \in {\rm Image} \ i_{- \alpha,
\varphi} = {\rm Kernel} \ \pi_{- \alpha, \varphi}$. Hence, $\xi :=
\pi_{- \alpha,\varphi} \circ \Phi (\gamma) = 0$. Assume by
contradiction that there exists $\delta$ satisfying
$\Pi(\delta,\gamma) \neq 0$ such that $c(\delta,\varphi^{-1}) <
\alpha$. We observe that $\Phi (\delta) \in {\rm Image} \ i_{
\alpha,\varphi^{-1}}$. Let $\eta \in
HF_{2n-*}^{(-\infty,\alpha]}(\varphi^{-1})$ be such that $\Phi (
\delta ) = \ i_{ \alpha,\varphi^{-1}} (\eta)$. It follows
from~$\ref{CDF3}$ that $ \Pi(\delta,\gamma) = L'( \eta , \xi) =
0 $. This contradicts the above assumption that
$\Pi(\delta,\gamma) \neq 0$. Thus we must have
$c(\delta,\varphi^{-1}) \geq \alpha$ for every $\delta$ satisfying
$\Pi(\delta,\gamma) \neq 0$. Hence,
$$ \inf_{_{\delta : \Pi(\delta,\gamma) \neq 0}} c(\delta,\varphi^{-1})
\ \geq \alpha = - \varepsilon - c(\gamma,\varphi).$$ Again, since
this inequality holds for every $\varepsilon > 0$ we conclude that
$$ \inf_{_{\delta : \Pi(\delta,\gamma) \neq 0}} c(\delta,\varphi^{-1})
\ \geq  - c(\gamma,\varphi).$$
\end{enumerate}
The proof is now complete.
\end{proof}

\bibliographystyle{gtart} \bibliography{link}

\end{document}